\documentclass[draft]{amsart}
\usepackage{amssymb,graphicx,enumerate,amsthm}
\usepackage{multirow}
\usepackage{xcolor}
\usepackage{cite}

\newcommand{\T}{\mathcal{T}}
\renewcommand{\S}{\mathcal{S}}
\newcommand{\Z}{\mathbb{Z}}
\newcommand{\F}{\mathcal{F}}
\newcommand{\G}{\mathcal{G}}

\numberwithin{equation}{section}
\newtheorem{theorem}{Theorem}[section]

\newtheorem{corollary}[theorem]{Corollary}

\begin{document}

\makeatletter
\def\imod#1{\allowbreak\mkern10mu({\operator@font mod}\,\,#1)}
\makeatother

\author{Alexander Berkovich}
   	\address{Department of Mathematics, University of Florida, 358 Little Hall, Gainesville FL 32611, USA}
   	\email{alexb@ufl.edu}

\author{Ali Kemal Uncu}
   \address{Research Institute for Symbolic Computation, Johannes Kepler University, Linz. Altenbergerstrasse 69
A-4040 Linz, Austria}
   \email{akuncu@risc.jku.at}

\thanks{Research of the first author is partly supported by the Simons Foundation, Award ID: 308929. Research of the second author is supported by the Austrian Science Fund FWF, SFB50-07 and SFB50-09 Projects.}

\title[\scalebox{.9}{Refined \lowercase{$q$}-Trinomial Coefficients and Two Infinite Hierarchies of \lowercase{$q$}-Series Identities}]{Refined \lowercase{$q$}-Trinomial Coefficients and Two Infinite Hierarchies of \lowercase{$q$}-Series Identities}
     
\begin{abstract} We will prove an identity involving refined $q$-trinomial coefficients. We then extend this identity to two infinite families of doubly bounded polynomial identities using transformation properties of the refined $q$-trinomials in an iterative fashion in the spirit of Bailey chains. One of these two hierarchies contains an identity which is equivalent to Capparelli's first Partition Theorem.
\end{abstract}

\keywords{Happy Birthday!, Capparelli's Partition Theorem, Infinite hierarchy, Refined $q$-trinomial coefficients; $q$-Series; Polynomial Identities}
  
\subjclass[2010]{11B65, 11C08, 11P81, 11P82, 11P83, 11P84,  05A10, 05A15, 05A17}


\date{\today}
   
\dedicatory{To Peter---king of the castle---Paule on the occasion of his 60th birthday}
   
\maketitle

\section{Introduction}

There are many important transformations for the $q$-binomial coefficients of the type \begin{equation}\label{Intro_identity}
\sum_{r=0}^L  \frac{q^{r^2}(q;q)_{2L}}{(q;q)_{L-r}(q;q)_{2r}} {2r \brack r-j}_q = q^{j^2} {2L \brack L-j}_q
\end{equation} in the $q$-series literature (see \cite{Berkovich_Warnaar}, and the references there). Throughout this work $|q|<1$. Here, the $q$-Pochhammer symbols are defined as\begin{align}
(a;q)_n &:= (1-a)(1-aq)(1-aq^2)\dots (1-aq^{n-1}),\\
\intertext{for any non-negative integer $n$. In addition, we have}
(a;q)_\infty&:= \lim_{n\rightarrow\infty}(a;q)_n,\\
(a_1,a_2,\dots,a_k;q)_n &:= (a_1;q)_n(a_2;q)_n\dots (a_k;q)_n.
\intertext{We extend the definition of $q$-Pochhamer symbols to negative $n$ using}
\label{Shifting_of_infinite_Products} (a;q)_n &= \frac{(a;q)_\infty}{(aq^n;q)_\infty}.\\
\intertext{Observe that \eqref{Shifting_of_infinite_Products} implies}
\label{One_over_neg_k} \frac{1}{(q;q)_n} &= 0\text{  if  }n <0.
\intertext{The $q$-binomial coefficients are defined as}
 \label{Binom_def}
  \displaystyle {m+n \brack m}_q &:= \left\lbrace \begin{array}{ll}\frac{(q;q)_{m+n}}{(q;q)_m(q;q)_{n}},&\text{for }m, n \geq 0,\\
   0,&\text{otherwise.}\end{array}\hspace{-.2cm}\right.
\end{align}

The salient features of \eqref{Intro_identity} are that the sum over a $q$-binomial coefficient multiplied by a simple factor yields a $q$-binomial coefficient, and the dependence on the variable $j$ is simple. Furthermore, this summation can be applied multiple times in an iterative fashion. This type of transformations were used by Bailey \cite{Bailey} and Slater \cite{Slater_List}, but the real value of the iterative power was first realized by Peter Paule \cite{Peter_Thesis, Peter85} and George E. Andrews \cite{Andrews_Multiple}.

For example, we start with the simple identity \begin{align}
\label{Kronecker_delta}  \delta_{L,0}&= \sum_{j=-L}^L (-1)^j q^{j\choose 2} {2L \brack L+j}_{q},
\end{align} where $\delta_{i,j}$ is the Kronecker delta function. Following the change of variable $L\mapsto r$ in \eqref{Kronecker_delta}, we multiply both sides by \begin{equation}
\label{what_to_multiply} \frac{q^{r^2}(q;q)_{2L}}{(q;q)_{L-r}(q;q)_{2r}},
\end{equation} and sum both sides with respect to $r\geq 0$. This yields \begin{equation}\label{First_Iteration}
 \frac{(q;q)_{2L}}{(q;q)_L}=\sum_{j=-L}^L (-1)^j q^{j^2+{j\choose 2}} {2L \brack L+j}_{q},
\end{equation} using the identity \eqref{Intro_identity}. 

It is well known that
\begin{align}
\label{Binom_limit}
\lim_{L\rightarrow\infty}{L\brack m}_q &= \frac{1}{(q;q)_m}.
\intertext{For any $j\in \mathbb{Z}_{\geq0}$ and $\nu =0$ or 1}
\label{Binom_limit2} \lim_{L\rightarrow\infty}{2L+\nu\brack L+j}_q &= \frac{1}{(q;q)_\infty}.
\end{align} 

Using \eqref{Binom_limit2} to take the limit $L\rightarrow\infty$ of \eqref{First_Iteration}, we get \begin{equation}
 (q;q)_{\infty}=\sum_{j=-\infty}^\infty (-1)^j q^{\frac{3j^2-j}{2}} .
\end{equation} This is nothing but Euler's Pentagonal Number Theorem. Note that this can also be viewed as a special case $(q,z)\mapsto(q^{3/2},-q^{1/2})$ of the celebrated result:
\begin{theorem}[Jacobi Triple Product Identity] For complex numbers $z\neq0$ and $|q|<1$, we have
\begin{equation}\label{JTP}
\sum_{j=-\infty}^\infty z^j q^{j^2} = \left(q^2,-zq,-\frac{q}{z};q^2\right)_\infty.
\end{equation} 
\end{theorem}

We can apply \eqref{Intro_identity} to the identity \eqref{First_Iteration} after changing the variable $L\mapsto r$ in \eqref{First_Iteration}, multiply both sides by \eqref{what_to_multiply} and sum both sides again with respect to $r\geq 0$. This yields \begin{equation}
\label{Second_iteration} \frac{(q;q)_{2L}}{(q;q)_L}\sum_{r= 0}^L q^{r^2} {L\brack r}_q=\sum_{j=-L}^L (-1)^j q^{2j^2+{j\choose 2}} {2L \brack L+j}_{q}.
\end{equation} Letting $L\rightarrow\infty$, using \eqref{Binom_limit}, \eqref{Binom_limit2}, \eqref{JTP} with $(q,z)\mapsto(q^{5/2}, -q^{1/2})$, and doing some simple simplifications one obtains the first Rogers--Ramanujan identity,
\[\sum_{r\geq 0} \frac{q^{r^2}}{(q;q)_r} = \frac{1}{(q,q^4;q^5)_\infty}.\]

Proceeding in the same fashion, one can keep on applying \eqref{Intro_identity} iteratively $\nu+1$ times to \eqref{Kronecker_delta}. This way one obtains the famous Andrews--Gordon infinite hierarchy of identities as $L\rightarrow\infty$,
\begin{equation}
\sum_{n_1,n_2,\dots ,n_\nu\geq0 } \frac{q^{N_1^2+N_2^2+\dots +N_\nu^2}}{(q;q)_{n_1}(q;q)_{n_2}\dots(q;q)_{n_\nu}} = \prod_{n \not\equiv 0, \pm (\nu+1)\text{ (mod $2\nu+3$)}} \frac{1}{1-q^n},
\end{equation} where $N_k = n_k + n_{k+1}+\dots +n_\nu\,$ for $k= 1,2,\dots ,\nu$.

In \cite{Andrews_Baxter}, Andrews and Baxter defined the $q$-trinomial coefficients,
\begin{align}
\label{Round_Tri_Def}
\left(\hspace{-.2cm}\begin{array}{c}L,\, b\\ a \end{array}\hspace{-.2cm};q\right)_2 &:= \sum_{n\geq 0 } q^{n(n+b)} \frac{(q;q)_L}{(q;q)_n(q;q)_{n+a}(q;q)_{L-2n-a}},\\
\label{T_n_Def} T_0\left( \hspace{-.2cm}\begin{array}{c}L\\a\end{array}\hspace{-.2cm};q\right) &:= q^{\frac{L^2-a^2}{2}} \left(\hspace{-.2cm}\begin{array}{c}L,\, a\\ a \end{array}\hspace{-.2cm};\frac{1}{q}\right)_2.
\end{align}
Following that Warnaar \cite{Warnaar_Refined, Warnaar_T} defined a refinement of these coefficients:
\begin{align}\label{T_Warnaar}\T&\left( \hspace{-.2cm}\begin{array}{c}L,\, M \\a,\, b\end{array}\hspace{-.2cm};q  \right) := \sum_{\substack{n\geq0, \\ L-a\equiv n \text{ (mod 2)}}} q^{\frac{n^2}{2}} {M\brack n}_q {M+b+\frac{L-a-n}{2} \brack M+b}_q {M-b +\frac{L+a-n}{2}\brack M-b}_q,\\
\label{S_Warnaar}\S&\left( \hspace{-.2cm}\begin{array}{c}L,\, M \\a,\, b\end{array}\hspace{-.2cm};q  \right) := \sum_{n\geq0} q^{n(n+a)} {M+L-a-2n\brack M}_q {M-a+b \brack n}_q {M+a-b \brack n+a}_q. \end{align}
These refined trinomials obey transformation properties somewhat similar to \eqref{Intro_identity}. Therefore, they can be used in an iterative fashion \cite{Warnaar_T}. 

In this paper, we prove a new doubly bounded polynomial identity using the symbolic tools developed by the Algorithmic Combinatorics group at the Research Institute for Symbolic Computation. 
\begin{theorem}\label{Initial_Seed_THM} For $L$ and $M$ being non-negative integers, we have
\begin{equation}\label{Seed_EQN}
\sum_{\substack{m\geq 0,\\ L\equiv m \text{ (mod 2)}}} q^{m^2}{3M \brack m}_{q^2} {2M + \frac{L-m}{2}\brack 2M}_{q^6} = \sum_{j=-\infty}^{\infty} q^{3j^2+2j} \T\left( \hspace{-.2cm}\begin{array}{c}L,\, M \\j,\, j\end{array}\hspace{-.2cm};q^6  \right).
\end{equation}
\end{theorem}

Then we use transformation properties for the refined trinomials defined in \eqref{T_Warnaar} and \eqref{S_Warnaar} to obtain two new infinite hierarchies of $q$-series identities. An unusual feature of these identities is the presence of various $q$-factorial bases such as in the following theorem with bases $q^2$, $q^3$, $q^6$ etc.

\begin{theorem}\label{End_of_T_hierarchy} Let $\nu$ be a positive integer, and let $N_k = n_k+n_{k+1}+\dots +n_{\nu}$, for $k=1,2,\dots, \nu$. Then,
\begin{align}
\nonumber\sum_{n_1,n_2,\dots,n_\nu\geq 0} &\frac{q^{3(N^2_1+N^2_2+\dots + N^2_\nu)} (-q;q^2)_{3n_\nu}}{(q^6;q^6)_{n_1}(q^6;q^6)_{n_2}\dots (q^6;q^6)_{n_{\nu-1}}(q^6;q^6)_{2n_\nu}}\\
\label{End_of_T_hierarchy_EQN} &\hspace{2cm}= \frac{(-q^3;q^3)_\infty}{(q^{12};q^{12})_\infty} (q^{6(\nu+1)},-q^{3\nu+1},-q^{3\nu+5};q^{6(\nu+1)})_\infty.
\end{align}
\end{theorem}

This paper is structured as follows. Section~\ref{Section_background} has a short list of known identities that will be needed later. Section~\ref{Section_Pf_Seed} is totally dedicated to the proof of Theorem~\ref{Initial_Seed_THM}. In Section~\ref{Section_Limits_transformations}, we discuss the asymptotics of $q$-trinomial coefficients, and present two transformation formulas of Warnaar for the refined $q$-trinomial coefficients \eqref{T_Warnaar} and \eqref{S_Warnaar}. We also discuss an analog of the Bailey Lemma (Theorem~\ref{Bailey_T_to_THM}). In Section~\ref{Section_First_Hierarchy}, we apply \eqref{Bailey_T_to_T} to \eqref{Seed_EQN} in an iterative fashion. This application yields a doubly bounded infinite hierarchy. The asymptotic analysis and the proof of Theorem~\ref{End_of_T_hierarchy} are also given in Section~\ref{Section_First_Hierarchy}. In Section~\ref{Section_Second_Hierarchy}, we use the second transformation \eqref{Bailey_T_to_S} of Theorem~\ref{Bailey_T_to_THM}, which yields another doubly bounded hierarchy of polynomial identities, and do its asymptotic analysis. In this way we see a connection with the Capparelli partition theorem. In Section~\ref{Sec_Outlook} we briefly discuss variants of Theorem~\ref{Initial_Seed_THM}.

\section{$q$-Binomial Theorem and Its Corollaries}\label{Section_background}

\begin{theorem}[q-Binomial Theorem]\label{qBin_THM} For variables $a,\, q,$ and $z$,
\begin{equation}\label{qBin_EQN}\sum_{n\geq 0} \frac{(a;q)_n}{(q;q)_n}t^n = \frac{(at;q)_\infty}{(t;q)_\infty}.\end{equation}
\end{theorem}

Note that by setting $(a,t)\mapsto(q^{-L},-zq^L)$  in \eqref{qBin_EQN}, and using \[{L \brack n}_q = \frac{(q^{-L};q)_n}{(q;q)_n}(-1)^n q^{Ln -{n\choose 2}},\] we derive \begin{equation}
\label{qExonential_sum} \sum_{n\geq 0} q^{n \choose 2}z^n{L \brack n}_q = (-z;q)_L.
\end{equation} We remark that \eqref{qExonential_sum} implies \begin{equation}\label{Sigma_splitof_qExponential_Sum}
\sum_{\substack{n\geq 0,\\n\equiv \sigma\text{ (mod 2)}}} q^{n \choose 2}z^n{L \brack n}_q = \frac{(-z;q)_L+(-1)^\sigma (z;q)_L}{2},
\end{equation} where $\sigma =0$ or 1.
 
Another important corollary of the $q$-binomial theorem (Theorem~\ref{qBin_THM}) is the polynomial analog of the identity \eqref{JTP} \cite[p.49, Ex.1]{Theory_of_Partitions}.
\begin{equation}\label{Finite_JTP}
\sum_{j=-M}^M q^{j^2} z^j {2M \brack M+j}_{q^2} = \left(-zq,-\frac{q}{z};q^2\right)_M.
\end{equation}

Note that \eqref{Kronecker_delta} is a special case of \eqref{Finite_JTP} with $(q,z)\mapsto({q}^{1/2},-{q}^{1/2})$. Another special case of \eqref{Finite_JTP} with $(q,z)\mapsto(q^3,q^2)$ is \begin{equation}\label{1_5_6_Fin_JTP} \sum_{j=-M}^M q^{3j^2+2j} {2M \brack M+j}_{q^6} = (-q,-q^5;q^6)_M.\end{equation}

\section{Proof of Theorem~\ref{Initial_Seed_THM}}\label{Section_Pf_Seed}

Let 
\begin{align}
\label{G_def}\G\left( L,M,k,q \right) &:= q^{(L-2k)^2} {3M \brack L-2k}_{q^2} {2M+k \brack k}_{q^6} \\ \intertext{and}
\label{F_def}\F \left( L,M,k,j,q \right) &:= q^{3j^2+2j+3(L-j-2k)^2} {M\brack L-j-2k}_{q^6}{M+j+k\brack k}_{q^6}{M+k\brack k+j}_{q^6}. 
\end{align} 
Note that \[ \sum_{k\geq 0} \G\left( L,M,k,q \right)\text{ and }\sum_{j,k\geq 0}\F \left( L,M,k,j,q \right)\] are the left-hand and right-hand sides of \eqref{Seed_EQN}, respectively.

The Mathematica packages \textit{Sigma} \cite{Sigma} and \textit{qMultiSum} \cite{qMultiSum} (both implemented by the Algorithmic Combinatorics group at the Research Institute for Symbolic Computation) are both capable of finding and automatically proving recurrences for these functions. Here we start with the recurrences that \textit{qMultiSum} finds for the summands \eqref{G_def} and \eqref{F_def}: 
\begin{align*}
&{q}^{9+18\,M} \left( 1-{q}^{12+6\,L+6\,M} \right) \G \left( L,M,k,q \right) -{q}^{4+12\,M} \left( 1+{q}^{2}+{q}^{4} \right)  \left( -1+{q}^{18+6\,L+12\,M} \right) \G \left( L+1,M,k,q \right)\\[-1.5ex]\\
 &+{q}^{1+6\,M} \left( 1-{q}^{24+6\,L+18\,M} \right)  \left( 1+{q}^{2}+{q}^{4} \right) \G \left( L+2,M,k,q \right)- \left( 1-{q}^{12+24\,M} \right) \G \left( L+3,M+1,k,q
 \right)\\[-1.5ex]\\
 & + \left( 1-{q}^{30+6\,L+24\,M} \right) \G \left( L+3,M,k,q \right)  +{q}^{6+12\,M} \left( 1+{q}^{6} \right)  \left( 1-{q}^{6+12\,M} \right) \G \left( L+1,M+1,k-1,q \right) \\[-1.5ex]\\
 &=0
\end{align*}
and
\begin{align*}
&{q}^{9+18\,M} \left( 1-{q}^{12+6\,L+6\,M} \right)\F \left( L,M,k-1,j-1,q \right) +{q}^{4+12\,M} \left( 1-{q}^{18+6\,L+12\,M} \right)\F \left( L+1,M,k-1,j,q \right)\\[-1.5ex]\\
&+{q}^{6+12\,M} \left( 1-{q}^{18+6\,L+12\,M} \right)\F \left( L+1,M,k-1,j-1,q \right) - \left( 1-{q}^{12+24\,M} \right)\F \left( L+3,M+1,k,j-1,q \right)\\[-1.5ex]\\
& +{q}^{3+6\,M} \left( 1-{q}^{24+6\,L+18\,M} \right)\F \left( L+2,M,k,j-1,q \right) +{q}^{5+6\,M} \left( 1-{q}^{24+6\,L+18\,M} \right)\F \left( L+2,M,k,j-2,q \right)\\[-1.5ex]\\
&+{q}^{1+6\,M} \left( 1-{q}^{24+6\,L+18\,M} \right)\F \left( L+2,M,k-1,j,q \right)+ \left( 1-{q}^{30+6\,L+24\,M} \right)\F \left( L+3,M,k,j-1,q \right) \\[-1.5ex]\\
&+ \left( 1+{q}^{6} \right) {q}^{6+12\,M} \left( 1-{q}^{6+12\,M} \right)\F \left( L+1,M+1,k-1,j-1,q \right)\\[-1.5ex]\\
&+{q}^{8+12\,M} \left( 1-{q}^{18+6\,L+12\,M} \right)\F \left( L+1,M,k,j-2,q \right) =0.
\end{align*}
Once summed over the variable $k$ for $\G(L,M,k,q)$, and varibles $k$ and $j$ for $\F(L,M,k,j,q)$, we see that they satisfy the same recurrence,
\begin{align}
\nonumber &{q}^{9+18\,M} \left( 1-{q}^{12+6\,L+6\,M} \right) \hat{S} \left( L,M,q \right) -{q}^{4+12\,M} \left( 1+{q}^{2}+{q}^{4} \right)  \left( -1+{q}^{18+6\,L+12\,M} \right) \hat{S} \left( L+1,M,q \right)\\[-1.5ex]\nonumber\\
\label{Sum_rec}&+{q}^{1+6\,M} \left( 1-{q}^{24+6\,L+18\,M} \right)  \left( 1+{q}^{2}+{q}^{4} \right) \hat{S} \left( L+2,M,q \right)+ \left( 1-{q}^{30+6\,L+24\,M} \right) \hat{S} \left( L+3,M,q \right) \\[-1.5ex]\nonumber\\
\nonumber &+{q}^{6+12\,M} \left( 1+{q}^{6} \right)  \left( 1-{q}^{6+12\,M} \right) \hat{S} \left( L+1,M+1,q \right) - \left( 1-{q}^{12+24\,M} \right) \hat{S} \left( L+3,M+1,q \right) =0.
\end{align}
This is also the same recurrence one would get from the package \textit{Sigma}. It remains to show that the the left-hand side and the right-hand side of \eqref{Seed_EQN} satisfy the same initial conditions.
Observe that \begin{equation}\label{Boundary_S_1}\hat{S}(L,0,q) = \frac{1+(-1)^L}{2},\text{  and  } \hat{S}(0,M,q)=1,\end{equation} for any non-negative integer $L$ and $M$.
Moreover, we have \begin{equation}\label{Boundary_S_2}\hat{S}(1,M,q) = q{3M\brack 1}_{q^2},\text{  and  } \hat{S}(2,M,q)={2M+1\brack 1}_{q^6}+q^4{3M \brack 2}_{q^2},\end{equation} for any non-negative integer $M$.
The recurrence \eqref{Sum_rec}, and the boundary conditions \eqref{Boundary_S_1} and \eqref{Boundary_S_2} uniquely define $\hat{S}(L,M,q)$. \hfill $\square$

\section{Asymptotics and Transformations of the Refined Trinomial Coefficients}\label{Section_Limits_transformations}

For $\sigma =0$ or 1, we have the following limits.
\begin{align}
\label{T_M_limit}\lim_{M\rightarrow\infty}\T\left( \hspace{-.2cm}\begin{array}{c}L,\, M \\a,\, b\end{array}\hspace{-.2cm};q  \right) &= \frac{1}{(q;q)_L}T_0 \left(\hspace{-.2cm}\begin{array}{c} L\\a \end{array}\hspace{-.2cm};q\right),\\
\label{T_L_limit}\lim_{\substack{L\rightarrow\infty,\\L-a\equiv \sigma\text{ (mod 2)}}}\T\left( \hspace{-.2cm}\begin{array}{c}L,\, M \\a,\, b\end{array}\hspace{-.2cm};q  \right) &= \frac{(-q^{1/2};q)_M + (-1)^\sigma (q^{1/2};q)_M}{2(q;q)_{2M}} {2M\brack M-b}_q,\\
\label{Tn_Trinomial_L_lim}\lim_{\substack{L\rightarrow\infty,\\ L-a\equiv \sigma\text{ (mod 2)}}} T_0 \left( \hspace{-.2cm}\begin{array}{c}L\\a\end{array}\hspace{-.2cm};q  \right) &= 
\frac{(-q^{1/2};q)_\infty + (-1)^\sigma (q^{1/2};q)_\infty}{2(q;q)_{\infty}}.\\
\intertext{Moreover,}
\label{S_M_limit}\lim_{M\rightarrow\infty}\S\left( \hspace{-.2cm}\begin{array}{c}L,\, M \\a,\, b\end{array}\hspace{-.2cm};q  \right) &=\frac{1}{(q;q)_L}\left(\hspace{-.2cm}\begin{array}{c} L,a\\a \end{array}\hspace{-.2cm};q\right)_2,\\
\label{S_L_limit}\lim_{L\rightarrow\infty}\S\left( \hspace{-.2cm}\begin{array}{c}L,\, M \\a,\, b\end{array}\hspace{-.2cm};q  \right) &=\frac{1}{(q;q)_{M}}{2M\brack M-b}_q,\\
\label{RoundTri_L_limit}\lim_{L\rightarrow\infty}\left(\hspace{-.2cm}\begin{array}{c} L,a\\a \end{array}\hspace{-.2cm};q\right)_2 &= \frac{1}{(q;q)_\infty}.
\end{align}

We would like to note that the limits \eqref{T_M_limit}, \eqref{S_M_limit}, and \eqref{S_L_limit} can be found in Warnaar's work \cite[(2.12),(2.13),(2.17)]{Warnaar_T}. The limit \eqref{Tn_Trinomial_L_lim} appears in Andrews--Baxter work \cite[(2.55),(2.56)]{Andrews_Baxter}. The limit \eqref{RoundTri_L_limit} is also discussed in \cite[(2.48)]{Andrews_Baxter}. The authors could not find the limit \eqref{T_L_limit} in the literature. This limit can be proven by using \eqref{Binom_limit} followed up with \eqref{Sigma_splitof_qExponential_Sum}.

Letting $M\rightarrow\infty$ in \eqref{Seed_EQN}, and using \eqref{T_M_limit}, we get
\begin{equation}\label{Thm_3.9_Cap2}
\sum_{\substack{m\geq0,\\ L\equiv m\text{ (mod 2)}}} q^{m^2} \frac{(q^6;q^6)_L}{(q^2;q^2)_m (q^6;q^6)_{(L-m)/2}} = \sum_{j=-L}^L q^{3j^2+2j} T_0\left( \hspace{-.2cm}\begin{array}{c} L\\j \end{array}\hspace{-.2cm} ; q^6\right).
\end{equation} Observe that after the change of variables $n=(L-m)/2$, this identity becomes \cite[(3.9)]{BerkovichUncu8} with $q\mapsto q^2$.

Replacing $L\mapsto 2L+\sigma$, with $\sigma = 0,1$, letting $L$ tend to $\infty$, we get the following with the aid of \eqref{T_L_limit},
\begin{align}
\nonumber(-q;q^2)_{3M} &+(-1)^\sigma (q;q^2)_{3M} \\ \label{T_limit_L_of_the_Seed}
&=(-q^{3};q^6)_M \sum_{j=-M}^M q^{3j^2+2j} {2M\brack M+j}_{q^6} + (-1)^\sigma (q^{3};q^6)_M \sum_{j=-M}^M (-1)^j q^{3j^2+2j} {2M\brack M+j}_{q^6}.
\end{align}

It is easy to check that \eqref{T_limit_L_of_the_Seed} follows from the identity \eqref{1_5_6_Fin_JTP}.

\begin{theorem}[Warnaar \cite{Warnaar_Refined,Warnaar_T}]For $L,M,a,b\in \Z$ and $ab\geq 0$
 \begin{align}
\label{T_to_T} \sum_{i=0}^M q^{\frac{i^2}{2}} {L+M-i\brack L}_q \T\left( \hspace{-.2cm}\begin{array}{c}L-i,\, i \\a,\, b\end{array}\hspace{-.2cm};q  \right) = q^{\frac{b^2}{2}} \T\left( \hspace{-.2cm}\begin{array}{c}L,\, M \\a+b,\, b\end{array}\hspace{-.2cm};q  \right).
\intertext{For $L,M,a,b\in \Z$ with $ab\geq 0$, and $|a|\leq M$ if $|b|\leq M$ and $|a+b|\leq L$, then}
\label{T_to_S}\sum_{i=0}^M q^{\frac{i^2}{2}} {L+M-i\brack L}_q \T\left( \hspace{-.2cm}\begin{array}{c}i,\, L-i \\b,\, a\end{array}\hspace{-.2cm};q  \right) = q^{\frac{b^2}{2}} \S\left( \hspace{-.2cm}\begin{array}{c}L,\, M \\a+b,\, b\end{array}\hspace{-.2cm};q  \right).
\end{align}\end{theorem}

The transformation formulas \eqref{T_to_T} and \eqref{T_to_S} directly imply the following theorem.

\begin{theorem}\label{Bailey_T_to_THM}Let $F_{L,M}(q)$ and $\alpha_j(q)$ be sequences, and $L,M,m,n \in \Z_{\geq 0}$. If 
\begin{align}
F_{L,M}(q) &= \sum_{j=-\infty}^\infty \alpha_j(q)\T\left( \hspace{-.2cm}\begin{array}{c}L,\, M \\mj,\, nj\end{array}\hspace{-.2cm};q  \right)\\
\intertext{holds, then}
\label{Bailey_T_to_T}\sum_{i\geq 0} q^{\frac{i^2}{2}}{L+M-i \brack L}_{q} F_{L-i,i}(q) &=  \sum_{j=-
\infty}^\infty q^{\frac{(nj)^2}{2}} \alpha_j(q)\T\left( \hspace{-.2cm}\begin{array}{c}L,\, M \\(m+n)j,\, nj\end{array}\hspace{-.2cm};q  \right)
\intertext{and}
\label{Bailey_T_to_S}\sum_{i\geq 0} q^{\frac{i^2}{2}}{L+M-i \brack L}_{q} F_{i,L-i}(q) &=  \sum_{j=-
\infty}^\infty q^{\frac{(mj)^2}{2}} \alpha_j(q)\S\left( \hspace{-.2cm}\begin{array}{c}L,\, M \\(m+n)j,\, mj\end{array}\hspace{-.2cm};q  \right)
\end{align}
are true.
\end{theorem}

Note that \eqref{Bailey_T_to_T} can be used in combination with an appropriately chosen identity in an iterative fashion. This leads to an infinite hierarchy of identities. On the other hand, the identity \eqref{Bailey_T_to_S} can only be used once. 

\section{The First Doubly Bounded Infinite Hierarchy and Its Asymptotics}\label{Section_First_Hierarchy}

We use \eqref{Bailey_T_to_T} $\nu$ times with $q\mapsto q^6$ on \eqref{Seed_EQN} and obtain the following infinite hierarchy.

\begin{theorem} Let $\nu$ be a positive integer, and let $N_k = n_k+n_{k+1}+\dots +n_{\nu}$, for $k=1,2,\dots, \nu$. Then,
\begin{align}
\label{T_hierarchy}\sum_{\substack{m,n_1,n_2,\dots,n_\nu\geq 0,\\ L+m\equiv N_1+N_2+\dots +N_\nu\text{ (mod 2)}}} &q^{m^2+3(N_1^2+N_2^2\dots+N_\nu^2)}{L+M-N_1\brack L}_{q^6} {3n_\nu \brack m}_{q^2}\\\nonumber
&\times{2n_\nu + \frac{L-m-N_1-N_2-\dots-N_\nu}{2}\brack 2n_\nu}_{q^6} \prod_{j=1}^{\nu-1}{L-\sum_{l=1}^{j} N_l+n_j \brack n_j}_{q^6}\\\nonumber
& = \sum_{j=-\infty}^\infty q^{3(\nu+1)j^2+2j} \T\left( \hspace{-.2cm}\begin{array}{c}L,\, M \\ (\nu+1)j,\, j\end{array}\hspace{-.2cm};q^6  \right).
\end{align}
\end{theorem}

We replace $L\mapsto 2L+\sigma$, with $\sigma=0,1$, and sum over $\sigma$, in \eqref{T_hierarchy}. Letting $L\rightarrow \infty$ and using the \eqref{qExonential_sum} and \eqref{T_L_limit}, we get the following theorem.
\begin{theorem}Let $\nu$ be a positive integer, and let $N_k = n_k+n_{k+1}+\dots +n_{\nu}$, for $k=1,2,\dots, \nu$. Then,
\begin{align}
\nonumber\sum_{n_1,n_2,\dots,n_\nu\geq 0} &\frac{q^{3(N_1^2 + N_2^2+\dots + N_\nu^2)}(-q;q^2)_{3n_\nu}}{(q^6;q^6)_{M-N_1}(q^6;q^6)_{n_1}(q^6;q^6)_{n_2}\dots (q^6;q^6)_{n_{\nu-1}}(q^6;q^6)_{2n_\nu}}\\
\label{T_hierarchy_Limit_L}&\hspace{2cm}= \frac{(-q^3;q^6)_M}{(q^6;q^6)_{2M}} \sum_{j=-M}^M q^{3(\nu+1)j^2 + 2j} {2M \brack M+j}_{q^6}.
\end{align}
\end{theorem}

The $\nu=1$ case of the identity \eqref{T_hierarchy_Limit_L} yields a finite analog of the identity \cite[(6.7)]{BerkovichUncu7}.
\begin{corollary}
\begin{equation}\label{Fin_Analog_of_6.7_Cap1}
\sum_{n=0}^M \frac{q^{3n^2}(-q;q^2)_{3n}}{(q^6;q^6)_{M-n}(q^6;q^6)_{2n}} =\frac{(-q^3;q^6)_M}{(q^6;q^6)_{2M}} \sum_{j=-M}^M q^{6j^2+2j} {2M\brack M+j}_{q^6}.
\end{equation}
\end{corollary}

We can also take the limit $M\rightarrow\infty$ in the identity \eqref{T_hierarchy}. Using \eqref{T_M_limit} we get another infinite hierarchy.

\begin{theorem}Let $\nu$ be a positive integer, and let $N_k = n_k+n_{k+1}+\dots +n_{\nu}$, for $k=1,2,\dots, \nu$. Then,
\begin{align}
\nonumber\sum_{\substack{m,n_1,n_2,\dots,n_\nu\geq 0,\\L+m\equiv N_1+N_2+\dots + N_\nu \text{ (mod 2)}}} q^{m^2+3(N_1^2+N_2^2+\dots +N_\nu^2)} &{3n_\nu \brack m}_{q^2}\\
\label{T_hierarchy_Limit_M}&\hspace{-2.5cm}\times{2n_\nu + \frac{L-m-N_1-N_2-\dots-N_\nu}{2}\brack 2n_\nu}_{q^6}\prod_{j=1}^{\nu-1} {L-\sum_{l=1}^j N_l +n_j \brack n_j}_{q^6}\\
\nonumber&\hspace{-2.5cm}= \sum_{j=-\infty}^{\infty} q^{3(\nu+1)j^2 +2j} T_0 \left(\hspace{-.2cm}\begin{array}{c}L\\ (\nu+1)j\end{array}\hspace{-.2cm};q^6 \right).
\end{align}
\end{theorem}

Letting $M\rightarrow\infty$ in \eqref{T_hierarchy_Limit_L} and using \eqref{Binom_limit2} and \eqref{JTP} proves Theorem~\ref{End_of_T_hierarchy}.

\section{The Second Doubly Bounded Infinite Hierarchy and Its Asymptotics}\label{Section_Second_Hierarchy}

Now we look at the implications of \eqref{Bailey_T_to_S}. First we apply this identity with $q\mapsto q^3$ to \eqref{Seed_EQN} with $q^2\mapsto q$:

\begin{theorem} For $L$ and $M$ non-negative integers, we have
\begin{equation}\label{Nu_equals_0}
\sum_{\substack{i,m\geq 0,\\ i+m\equiv 0\text{ (mod 2)}}} q^{\frac{m^2+3i^2}{2}} {L+M-i\brack L}_{q^3}{3(L-i)\brack m}_{q}{2(L-i)+\frac{i-m}{2}\brack 2(L-i)}_{q^3}=\sum_{j=-\infty}^\infty q^{3j^2+j} \S\left( \hspace{-.2cm}\begin{array}{c}L,\, M \\ 2j,\, j\end{array}\hspace{-.2cm};q^3  \right).
\end{equation}
\end{theorem}

We can also apply \eqref{Bailey_T_to_S} with $q\mapsto q^3$ to \eqref{T_hierarchy} with $q^2\mapsto q$. This yields the following result.

\begin{theorem}
 Let $\nu$ be a positive integer, and let $N_k = n_k+n_{k+1}+\dots +n_{\nu}$, for $k=1,2,\dots, \nu$. Then,
\begin{align}
\label{S_hierarchy}\sum_{\substack{i,m,n_1,n_2,\dots,n_\nu\geq 0,\\ i+m\equiv N_1+N_2+\dots +N_\nu\text{ (mod 2)}}} &q^{\frac{m^2+3(i^2+N_1^2+N_2^2\dots+N_\nu^2)}{2}}{L+M-i\brack L}_{q^3}{L-N_1 \brack i}_{q^3} {3n_\nu \brack m}_{q}\\\nonumber
&\times{2n_\nu + \frac{i-m-N_1-N_2-\dots-N_\nu}{2}\brack 2n_\nu}_{q^3} \prod_{j=1}^{\nu-1}{i-\sum_{l=1}^{j} N_l+n_j \brack n_j}_{q^3}\\\nonumber
& = \sum_{j=-\infty}^\infty q^{3{\nu+2 \choose 2}j^2+j} \S\left( \hspace{-.2cm}\begin{array}{c}L,\, M \\ (\nu+2)j,\, (\nu+1)j\end{array}\hspace{-.2cm};q^3  \right).
\end{align}
\end{theorem}

Taking the limit $M\rightarrow\infty$ in \eqref{Nu_equals_0}, and changing the summation variable $(i-m)/2=n$ we get 
\begin{equation}\label{Limit_Andres_K}
\sum_{n,m\geq 0} q^{Q(m,n)} {3(L-2n-m)\brack m}_{q}{2(L-2n-m)+n\brack n}_{q^3}=\sum_{j=-\infty}^\infty q^{3j^2+j} \left( \hspace{-.2cm}\begin{array}{c}L,\, 2j\\ 2j\end{array}\hspace{-.2cm};q^3  \right)_2,
\end{equation} where $Q(m,n) := 2m^2+6mn+6n^2$.

The polynomials on the right-hand side of \eqref{Limit_Andres_K} were first discussed by Andrews in \cite{Andrews_P_Capparelli}. The identity, on the other hand, was first proven in \cite{BerkovichUncu7}.

The limit $L\rightarrow \infty$ in \eqref{Nu_equals_0} yields
\begin{equation}\label{Limit_Andrews_DK}
\sum_{n,m\geq 0} q^{Q(m,n)} \frac{(q^3;q^3)_M}{(q;q)_m(q^3;q^3)_{n}(q^3;q^3)_{M-2n-m}}=\sum_{j=-\infty}^\infty q^{3j^2+j} {2M\brack M+j}_{q^3}.
\end{equation}

This formula was first discussed in \cite{BerkovichUncu7}, and it is proven in a wider context in \cite{BerkovichUncu8}.

Finally, when $L$ and $M$ both tend to $\infty$, a simple change of variables together with the Jacobi Triple Product identity \eqref{JTP} yields
\begin{equation}\label{KR1}
\sum_{m,n\geq 0} \frac{q^{Q(m,n)}}{(q;q)_m (q^3;q^3)_n} =  (-q^2,-q^4;q^6)_\infty(-q^3;q^3)_\infty,
\end{equation}  where $Q(m,n) =2m^2 + 6mn +6n^2$, after simplifications.

The identity \eqref{KR1} was recently proposed independently by Kanade--Russell \cite{Kanade_Russell} and Kur\c{s}ung\"oz \cite{Kagan}.  They showed that \eqref{KR1} is equivalent to the following partition theorem:

\begin{theorem}[Capparelli's First Partition Theorem \cite{Capparelli_proof}]\label{m1_Capparelli} 
For any integer $n$, the number of partitions of $n$ into distinct parts where no part is congruent to $\pm 1$ modulo $6$ is equal to the number of partitions of $n$ into parts, not equal to $1$, where the minimal difference between consecutive parts is 2; the difference between consecutive parts is greater than or equal to $4$ unless consecutive parts are $3k$ and $3k+3$ (yielding a difference of 3), or $3k-1$ and $3k+1$ (yielding a difference of 2) for some $k\in \mathbb{N}$.
\end{theorem}

Taking limits $M\rightarrow \infty$ and $L\rightarrow\infty$ in \eqref{S_hierarchy}, we get Theorems~\ref{S_hierarchy_lim_M} and \ref{S_hierarchy_lim_L}, respectively.

\begin{theorem}\label{S_hierarchy_lim_M} Let $\nu$ be a positive integer, and let $N_k = n_k+n_{k+1}+\dots +n_{\nu}$, for $k=1,2,\dots, \nu$. Then,
\begin{align}
\nonumber\sum_{\substack{i,m,n_1,n_2,\dots,n_\nu\geq 0,\\ i+m \equiv N_1+N_2+\dots + N_\nu \text{ (mod 2)}}} &q^{\frac{m^2+3(i^2+N_1^2+N_2^2+\dots + N_\nu^2)}{2}} {L-N_1\brack i}_{q^3}{3n_\nu\brack m}_{q}\\
\label{S_hierarchy_lim_M_EQN}&\hspace{-1cm}\times {2n_\nu + \frac{i-N_1-N_2-\dots-N_\nu-m}{2}\brack 2n_\nu}_{q^3} \prod_{j=1}^{\nu-1} {i- \sum_{k=1}^j N_k+n_j \brack n_j}_{q^3}\\
\nonumber&\hspace{-1cm}= \sum_{j=\infty} ^\infty q^{3{\nu+2\choose 2}j^2+j} \left(\hspace{-.2cm}\begin{array}{c} L,\, (\nu+2)j\\ (\nu+2)j \end{array}\hspace{-.2cm}; q^3 \right)_2.
\end{align}
\end{theorem}

\begin{theorem}\label{S_hierarchy_lim_L} Let $\nu$ be a positive integer, and let $N_k = n_k+n_{k+1}+\dots +n_{\nu}$, for $k=1,2,\dots, \nu$. Then,
\begin{align}
\nonumber\sum_{\substack{i,m,n_1,n_2,\dots,n_\nu\geq 0,\\ i+m \equiv N_1+N_2+\dots + N_\nu \text{ (mod 2)}}} &q^{\frac{m^2+3(i^2+N_1^2+N_2^2+\dots + N_\nu^2)}{2}} {M \brack i}_{q^3}{3n_\nu\brack m}_{q}\\
\label{S_hierarchy_lim_L_EQN}&\hspace{-1cm}\times {2n_\nu + \frac{i-N_1-N_2-\dots-N_\nu-m}{2}\brack 2n_\nu}_{q^3} \prod_{j=1}^{\nu-1} {i- \sum_{k=1}^j N_k+n_j \brack n_j}_{q^3}\\
\nonumber&\hspace{-1cm}= \sum_{j=\infty} ^\infty q^{3{\nu+2\choose 2}j^2+j} {2M\brack M+(\nu+1)j}_{q^3}.
\end{align}
\end{theorem}

Finally, by letting $L\rightarrow \infty$ in \eqref{S_hierarchy_lim_M_EQN}, and using \eqref{RoundTri_L_limit} and  \eqref{JTP}, we get the following result.

\begin{theorem}\label{End_of_s_hierarchy}  Let $\nu$ be a positive integer, and let $N_k = n_k+n_{k+1}+\dots +n_{\nu}$, for $k=1,2,\dots, \nu$. Then,
\begin{align}
\nonumber\sum_{\substack{i,m,n_1,n_2,\dots,n_\nu\geq 0,\\ i+m \equiv N_1+N_2+\dots + N_\nu \text{ (mod 2)}}} &\frac{q^{\frac{m^2+3(i^2+N_1^2+N_2^2+\dots + N_\nu^2)}{2}}}{(q^3;q^3)_i} {3n_\nu\brack m}_{q}\\
\label{End_of_T_hierarchy_EQN}&\hspace{-1cm}\times {2n_\nu + \frac{i-N_1-N_2-\dots-N_\nu-m}{2}\brack 2n_\nu}_{q^3} \prod_{j=1}^{\nu-1} {i- \sum_{k=1}^j N_k+n_j \brack n_j}_{q^3}\\\nonumber
 &\hspace{2cm}= \frac{(q^{6{\nu+2\choose 2}},-q^{3{\nu+2\choose 2}+1},-q^{3{\nu+2\choose 2}-1};q^{6{\nu+2\choose 2}})_\infty}{(q^{3};q^{3})_\infty}.
\end{align}
\end{theorem}

Note that Theorem~\ref{End_of_s_hierarchy} can also be proven by taking the limit $M\rightarrow\infty$ in \eqref{S_hierarchy_lim_L_EQN}, and using \eqref{Binom_limit2} together with \eqref{JTP}.

\section{Outlook}\label{Sec_Outlook}

We would like to note that the identity \eqref{Seed_EQN} is not an isolated incident. This shows that there is a more complex structure behind and there is much more to discover. We would like to give two such example identities that we prove similarly to Theorem~\ref{Initial_Seed_THM}. Let \begin{align*}
\T_{1}&\left( \hspace{-.2cm}\begin{array}{c}L,\, M \\a,\, b\end{array}\hspace{-.2cm};q  \right) :=  \sum_{\substack{n\geq0, \\ L-a\equiv n \text{ (mod 2)}}} q^{n\choose 2} {M\brack n}_q {M+b+\frac{L-a-n}{2} \brack M+b}_q {M-b +\frac{L+a-n}{2}\brack M-b}_q,\intertext{and}
\T_{-1}&\left( \hspace{-.2cm}\begin{array}{c}L,\, M \\a,\, b\end{array}\hspace{-.2cm};q  \right) :=  \sum_{\substack{n\geq0, \\ L-a\equiv n \text{ (mod 2)}}} q^{n+1\choose 2} {M\brack n}_q {M+b+\frac{L-a-n}{2} \brack M+b}_q {M-b +\frac{L+a-n}{2}\brack M-b}_q.
\end{align*}
Then we have the following theorem.

\begin{theorem}\label{Outro_Seeds_THM} For $L$ and $M$ being non-negative integers, we have
\begin{align}\label{Seed_EQNS}
\sum_{\substack{m\geq 0\\L\equiv m\ \text{(mod }2\text{)}}} q^{m^2 \mp m}{3M \brack m}_{q^2} {2M + \frac{L-m}{2}\brack 2M}_{q^6} &= \sum_{j=-\infty}^{\infty} q^{3j^2+j} \T_{\pm 1}\left( \hspace{-.2cm}\begin{array}{c}L,\, M \\j,\, j\end{array}\hspace{-.2cm};q^6  \right).
\end{align}
\end{theorem}

We are planning to address Theorem~\ref{Outro_Seeds_THM} and its implications elsewhere.

\section{Acknowledgement}

We would like to thank Research Institute for Symbolic Computation for the warm hospitality.

Research of the first author is partly supported by the Simons Foundation, Award ID: 308929. Research of the second author is supported by the Austrian Science Fund FWF, SFB50-07 and SFB50-09 Projects.

We thank Chris Jennings-Shaffer for his careful reading of the manuscript.

\end{document}